\RequirePackage{fix-cm}
\documentclass[12pt]{elsarticle} 

\usepackage{mathptmx}      

\usepackage[utf8]{inputenc}
\usepackage{epstopdf}

\usepackage{array}
\usepackage{amsfonts}
\usepackage{amsmath}
\usepackage{amsthm}
\usepackage{amssymb}
\usepackage{bm}
\usepackage{enumerate}
\usepackage{comment}

\usepackage{epsfig} 

\newcommand{\myvec}[1]{{\bm {#1}}}  
\newcommand{\mymat}[1]{{\bm {#1}}}  
\newcommand{\myqed}{} 

\newtheorem{lemma}{Lemma} 
\newtheorem{theorem}{Theorem} 
\theoremstyle{definition} 
\newtheorem{definition}{Definition} 

\journal{arXiv} 

\begin{document}

\begin{frontmatter} 

%

\title{On the Transition of Charlier Polynomials to the
Hermite Function}


\author{Martin Nilsson}

\address{
SICS (Swedish Institute of Computer Science),
POB 1263,
164 29 Kista, Sweden}
\ead{convergence.theorem@drnil.com}



\begin{abstract}
       It has been known for over 70 years that there is an asymptotic transition of Charlier polynomials
to Hermite polynomials. This transition, which is still presented in its classical form in modern reference works,
is valid if and {\em only} if a certain parameter is integer.
In this light, it is surprising that a much more powerful transition exists from
Charlier polynomials to the Hermite {\em function}, valid for any real value of the parameter.
This greatly strengthens the asymptotic connections between Charlier polynomials and special functions,
with applications for instance in queueing theory.

It is shown in this paper that the
convergence is uniform in bounded intervals, and a sharp rate bound is proved.
It is also shown that there is a transition of derivatives of Charlier polynomials to the derivative of
the Hermite function, again with a sharp rate bound. Finally, it is proved that
zeros of Charlier polynomials converge to zeros of the Hermite function.
While rigorous, the proofs use only elementary techniques.


\end{abstract}


\end{frontmatter}    




\section{Introduction}

A unique feature of Charlier polynomials
\cite{
Charlier.1905udf,
erdelyi.et.al.1953htf,
Szego.1975op,
olver.et.al.2010nho,
Koekoek.et.al.2010hop}
is their affinity with the Poisson distribution. This has many important applications.
Charlier polynomials concisely express the behaviour of
{\em Erlang loss systems}, a fundamental concept in queueing theory
\cite{ Jagerman.1974spo, karlin.mcgregor.1958msq,Kijima.1990otl}.
Another example is the generalization of stochastic integrals over Poisson
distributions to multiple stochastic integrals,
which can be effectively computed using Charlier polynomials
\cite{Engel.1982tms,
Xiu.2010nmf},
while a third example is that of random matrices over Poisson distributions
\cite{konig.2005ope},
which can be characterized by Charlier polynomial zeros.

\begin{figure}
\includegraphics[clip, trim=1.0cm 0.5cm 2.3cm 15cm ]{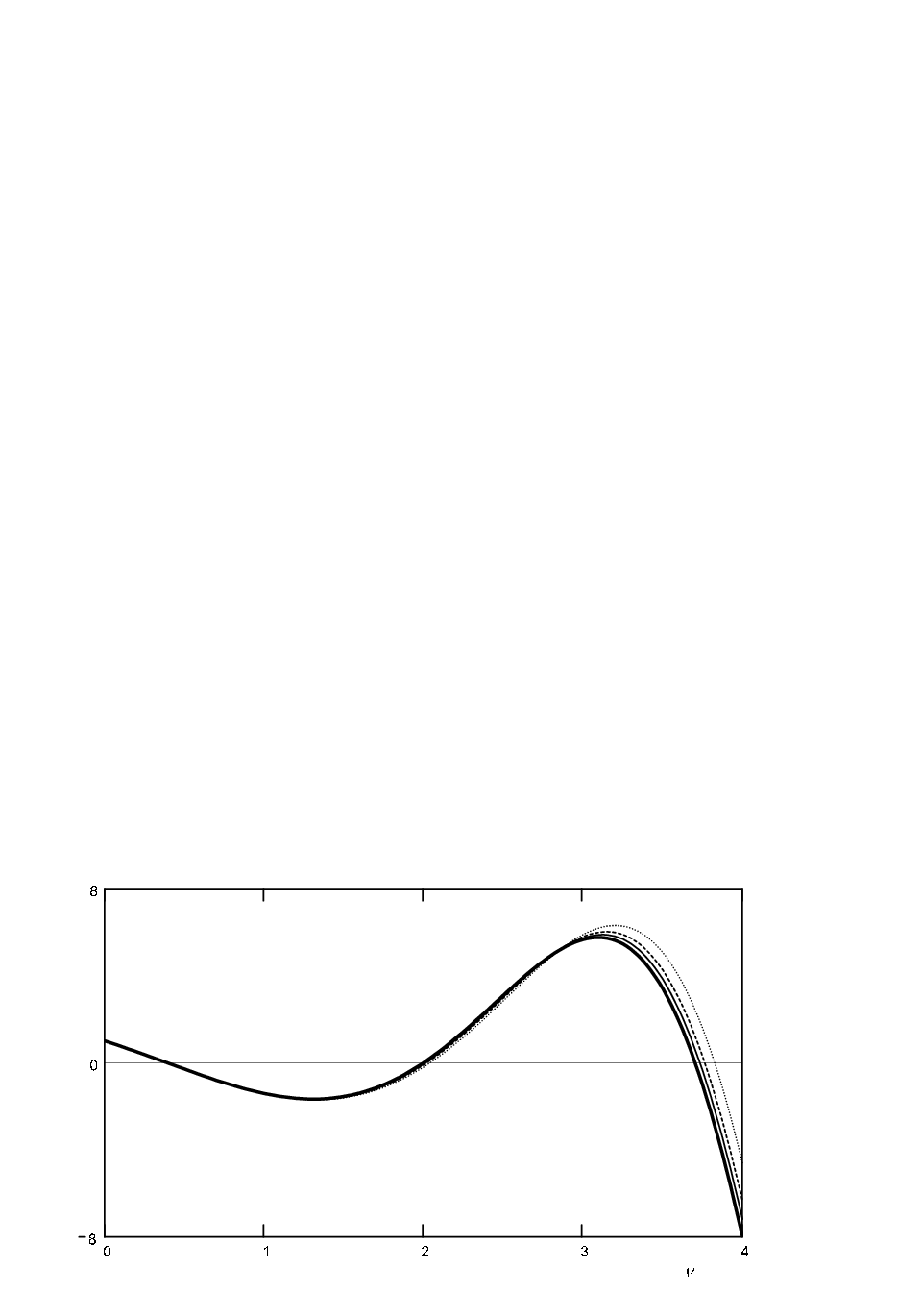}
\caption{Transition of Charlier polynomials 
${{\left( 2a \right)}^{\nu/2}}C_{\lceil a+\sqrt{a}\rceil}(\nu,a)$ 
(thin lines) 
to the Hermite function $H_\nu(-1/\sqrt{2})$ (thick line). The different values of the parameter $a$ are 100 (dotted line), 400 (dashed line), and 1600 (solid thin line).}
\label{fig:1}       
\end{figure}

High-dimensional or asymptotic problems typically engage Charlier polynomials of high degree and order (index).
For instance, the asymptotic behavior in the number of servers of Erlang loss systems is described by Charlier polynomials
whose degree and order tend to infinity simultaneously according to the {\em Halfin-Whitt regime} \cite{Halfin.Whitt.1981htl}.
At a first glance, the classical formula
(\cite[Eq. 9.14.12]{Koekoek.et.al.2010hop},
\cite[p. 532]{meixner.1939efd},
\cite[Eq. 18.21.9]{olver.et.al.2010nho},
\cite[Eq. 2.82.7]{Szego.1975op})
\begin{equation}\label{usual_expression_of_old_equation}
\underset{a\to \infty }{\mathop{\lim }}\,{{\left( 2a \right)}^{n/2}}C_{n}(a+x\sqrt{2a},a)={{\left( -1 
\right)}^{n}}{{H}_{n}}(x)
\end{equation}
appears useful for reducing Charlier polynomials in this limit,
but unfortunately, this formula holds {\em only} for non-negative integer $n$.
In light of the long standing of this formula, it can somewhat surprisingly be shown that
\begin{equation}\label{old_identity} 
\underset{a\to \infty }{\mathop{\lim }}\,{{\left( 2a \right)}^{\nu/2}}C_{\lceil a-
x\sqrt{2a}\rceil}(\nu,a)=H_\nu(x)
\end{equation}
for {\em any real $x$ and $\nu$}, a much stronger statement (fig.\ref{fig:1}).
Here, the ceiling function $\lceil x \rceil$ denotes the smallest integer not smaller 
than $x$, and $H_\nu(x)$ denotes the Hermite function~\cite[Ch.~10]{lebedev.1972sfa}.

A proof of~\eqref{usual_expression_of_old_equation} has been given via Krawtchouk
polynomials~\cite[pp.~36-37]{Szego.1975op}. When $\nu$ is non-negative real and $a-x\sqrt{2a}$ is integer,
pointwise convergence of \eqref{old_identity}
(without rate bound) follows implicitly from~\cite{dominici.2007aao}.

In section \ref{section:conv}, it is proved rigorously that convergence to the Hermite function holds for any 
real $\nu$, and that convergence is uniform for $\nu$ and $x$ in any 
bounded interval. A sharp rate bound is established.
The same technique is then employed in section
\ref{section:derivative} to prove 
that there is a similar transition of the derivative with respect to $\nu$.
Also here, a sharp rate bound is provided. These results are used in
section \ref{section:zeros} for proving that zeros of Charlier polynomials converge
to zeros of the Hermite function.

Below is first a recollection of some well-known definitions and recurrence relations from
\cite{
erdelyi.et.al.1953htf,
Koekoek.et.al.2010hop,
lebedev.1972sfa,
olver.et.al.2010nho}
in order to make the paper self-contained. This is followed by three sections, each proving
an aspect of the transition of Charlier polynomials to Hermite functions.

The notation ``$A\triangleq B$'' is used for ``$A$ is defined $B$'', in order to make the 
introduction of new symbols more explicit. The expression ``bounded $\nu \le -3$'' is 
shorthand for ``$\nu$ in any bounded interval $[{{\nu }_{0}},-3]$''.
We will abbreviate Charlier polynomials $C_n(x,a)$ as $c_{n}^{a}(x)$,
$c_{n}(x)$, $c_{n}^{a}$, or even  $c_{n}$, unless there is a risk for misunderstanding.
They
can be defined for positive $a$ and non-negative integer $n$ by
\cite[Eq. 10.25.4]{erdelyi.et.al.1953htf},
\cite[Eq. 18.20.8]{olver.et.al.2010nho},
\begin{equation}\label{def_charlier_polynomial}
   c_{n}^{a}\left( x \right)\triangleq \underset{k=0}{\overset{n}{\mathop 
\sum }}\,\left( \begin{matrix}
   n  \\
   k  \\
\end{matrix} \right)\left( \begin{matrix}
   x  \\
   k  \\
\end{matrix} \right)k!{{\left( -a \right)}^{-k}}
\end{equation}
where 
\[
\left(
\begin{matrix}
   x  \\
   k  \\
\end{matrix} \right)\triangleq
\begin{cases}
   {x\left( x-1 \right)\cdot \ldots ~\cdot (x-k+1)}/{k!} & \text{for}~k\ge 1  \\
   1 & \text{for}~k=0  \\
\end{cases}
\]
These polynomials obey the three-term recurrence relation
\cite[Eq. 10.25.8]{erdelyi.et.al.1953htf},
\cite[Eq. 18.22.2]{olver.et.al.2010nho},
\begin{equation}\label{three_term_recurrence_rule} 
-xc_{n}^{a}\left( x \right)=ac_{n+1}^{a}\left( x \right)-\left( n+a \right)c_{n}^{a}\left( x 
\right)+nc_{n-1}^{a}(x)
\end{equation}
and the difference equation
\cite[Eq. 10.25.9]{erdelyi.et.al.1953htf},
\cite[Eq. 18.22.12]{olver.et.al.2010nho},
\begin{equation}\label{difference_equation} 
-nc_{n}^{a}\left( x \right)=ac_{n}^{a}\left( x+1 \right)-\left( x+a \right)c_{n}^{a}\left( x 
\right)+xc_{n}^{a}(x-1)
\end{equation}
as well as the backward recurrence relation
\cite[Eq. 9.14.8]{Koekoek.et.al.2010hop},
\begin{equation}\label{backward_recurrence_relation} 
\frac{x}{a}c_{n-1}^{a}\left( x-1 \right)=c_{n-1}^{a}(x)-c_{n}^{a}(x)
\end{equation}
The Hermite function ${{H}_{\nu }}(x)$ is a solution of the differential equation~\cite[Eq. 10.2.3]{lebedev.1972sfa} 
\begin{equation}\label{Hermite_equation} 
   {y}''=2x{y}'-2\nu y
\end{equation}
and satisfies the three-term recurrence~\cite[Eq. 10.4.7]{lebedev.1972sfa}
\begin{equation}\label{Hermite_three_term_recurrence} 
{{H}_{\nu +1}}\left( x \right)-2x{{H}_{\nu }}\left( x \right)+2\nu {{H}_{\nu -1}}\left( x 
\right)=0
\end{equation}
and the derivative rule~\cite[Eq. 10.4.4]{lebedev.1972sfa}
\begin{equation}\label{derivative_rule} 
{H'_{\nu }}\left( x \right)=2\nu {{H}_{\nu -1}}(x)
\end{equation}
It can be defined by~\cite[Eq. 10.2.8]{lebedev.1972sfa}
\begin{equation}\label{Hermite_function} 
{{H}_{\nu }}\left( x \right)\triangleq 
{{2}^{\nu }}\sqrt{\pi }\left[ \left( \frac{1}{\Gamma \left( \frac{1-
\nu }{2} \right)} \right)M\left( -\frac{\nu }{2};\frac{1}{2};{{x}^{2}} \right)-
\frac{2x}{\Gamma \left( -\frac{\nu }{2} \right)}M\left( \frac{1-
\nu }{2};\frac{3}{2};{{x}^{2}} \right) \right]
\end{equation}
where $M$ is the confluent hypergeometric function of the first kind. When the expression 
involves a gamma function of a non-positive integer argument, the expression should be 
interpreted by its limiting value.

\section{Transition of Charlier polynomials}
\label{section:conv}

\begin{theorem}
\label{theorem:1}
For real $x$, $\nu$, and positive $a$,
\[
{{\left( 2a \right)}^{\nu /2}}c_{\lceil a-x\sqrt{2a}\rceil}^{a}(\nu )={{H}_{\nu }}\left( x 
\right)+O\left( \frac{1}{\sqrt{a}} \right)
\]
where $c_{n}^{a}(\nu )$ are Charlier polynomials and ${{H}_{\nu }}(x)$ is the Hermite 
function. The error bound $O\left( {1}/{\sqrt{a}}\; \right)$ is uniform for $\nu$ and $x$ in 
any bounded interval, and is sharp in the sense that there are $\nu$ and $x$ such that the 
error is proportional to $1/\sqrt{a}$ for arbitrarily large $a$.
\end{theorem}

Proving asymptotic properties of Charlier polynomials is difficult, since these do not
satisfy a second-order linear ordinary differential equation with respect to the independent 
variable~\cite{dunster.2001uae}. However, the three-term recurrence relation~\eqref{three_term_recurrence_rule} is a discretization of such a 
differential equation~\eqref{Hermite_equation}. This can be used in order to prove the theorem in the following way: It 
is first proven for the special case $x=0$ and $\nu \le -4$ (Lemma~\ref{factor_p}-\ref{Rhead}), and then generalized 
to arbitrary real $\nu$ (Lemma~\ref{y_zero_case}). After that, the scaled polynomials are shown to 
approximate a Cauchy polygon converging to the ${{H}_{\nu }}\left( x \right)$ solution of the 
Hermite differential equation initial value problem (Lemma~\ref{Cauchy_pol_conv}).

\subsection{Convergence for $x=0$ and $\nu \le -4$}

For notational convenience, define $A\triangleq \lceil a \rceil$ and 
\begin{equation}\label{def_y} 
\begin{matrix}
   y_{\nu }^{a}(x)\triangleq {{\left( 2a \right)}^{\nu /2}}c_{\lceil a-x\sqrt{2a}\rceil}^{a}(\nu )  \\
\end{matrix}
\end{equation}
The superscript will be left out in $y_{\nu }^{a}$ and $c_{n}^{a}$ unless there is a risk 
for misunderstanding. Consider the case $x=0$ and $\nu \le -4$. By the definition of Charlier 
polynomials~\eqref{def_charlier_polynomial}, 
\[
   {{c}_{n}}\left( \nu  \right)=\underset{k=0}{\overset{n}{\mathop \sum }}\,\left( \begin{matrix}
   n  \\
   k  \\
\end{matrix} \right)\frac{\Gamma \left( k-\nu  
\right)}{\Gamma \left( -\nu  \right)}~{{a}^{-k}}
\]
In order to prove that
	\[\underset{a\to \infty }{\mathop{\lim }}\,{{y}_{\nu }}(0)={{H}_{\nu }}\left( 0 
\right)=\frac{{{2}^{\nu }}\sqrt{\pi }}{\Gamma \left( \frac{1-\nu }{2} \right)}\]
using the definition of ${{H}_{\nu }}(x)$ in~\eqref{Hermite_function}, ${{y}_{\nu }}(0)$ can be expressed as a 
sum
\begin{equation}\label{y_sum} 
{{y}_{\nu }}(0)=\frac{{{2}^{\nu /2}}}{\Gamma \left( -\nu  
\right)}\underset{k=0}{\overset{A}{\mathop \sum }}\,{{T}_{k}}
\end{equation}
of terms
\[
{{T}_{k}}\triangleq {{a}^{\nu /2}}\frac{\Gamma \left( k-\nu  
\right)}{k!}~\frac{A!{{a}^{-k}}}{\left( A-k \right)!}
\]
When $\nu$ is negative, these are all positive. The series is difficult to sum due to multiple levels of
numerical cancellation, but can be estimated by separating the factors. Another difficulty is the changing behaviour of $T_k$
with increasing $a$. This problem can be remedied by defining a border between "head" and "tail" sections that
increases with a properly tuned power of $a$.

\begin{lemma}\label{factor_p} 
The factor
\[
p\left( k \right)\triangleq \frac{A!{{a}^{-k}}}{\left( A-k \right)!}={{\left( \frac{a}{A} \right)}^{-
k}}~\underset{j=0}{\overset{k-1}{\mathop \prod }}\,\left( 1-\frac{j}{A} \right)
\]
for $1\le k\le A~$satisfies
\[
p\left( k \right)\le \exp\left( -\frac{{{k}^{2}}}{2A} \right)\left[ 1+O\left( \frac{k}{a} \right) 
\right]
\]
and for $1\le k<A/2$,
\[
   p\left( k \right)\ge \exp\left( -\frac{{{k}^{2}}}{2A} 
\right)\left[ 1+O\left( \frac{k}{a}+\frac{{{k}^{3}}}{{{a}^{2}}} \right) \right]
\]
\end{lemma}

\begin{proof} 
Define the ``nuisance factor'' due to truncation by the ceiling function by
\begin{equation}\label{nuisance_factor} 
\beta \triangleq {{\left( \frac{a}{A} \right)}^{k}}={{\left( \frac{a}{a+\left( \lceil a\rceil-a \right)} 
\right)}^{k}}={{\left( \frac{a}{a+\theta } \right)}^{k}}=1+O\left( \frac{k}{a} \right)
\end{equation}
where $0\le \theta <1$. For $1\le k\le A$, taking the logarithm of $\beta p(k)$ and Taylor 
expanding,
\begin{align}
 \ln \left[\beta p\left( k \right)\right] & = \underset{j=0}{\overset{k-1}{\mathop \sum }}\,\ln \left( 1-
\frac{j}{A} \right) \notag\\ 
 & = \underset{j=0}{\overset{k-1}{\mathop \sum }}\,\left( -\frac{j}{A}-
\frac{{{j}^{2}}}{2{{A}^{2}}}-\frac{{{j}^{3}}}{3{{A}^{3}}}-\frac{{{j}^{4}}}{4{{A}^{4}}}-\cdots  \right) 
\notag\\ 
 & = -\frac{k\left( k-1 \right)}{2A}-\underset{j=0}{\overset{k-1}{\mathop 
\sum }}\,\left( \frac{{{j}^{2}}}{2{{A}^{2}}}+\frac{{{j}^{3}}}{3{{A}^{3}}}+\frac{{{j}^{4}}}{4{{A}^{4}}
}+\cdots  \right) \notag\\ 
 & \triangleq  -\frac{k\left( k-1 \right)}{2A}-{{R}_{p} } \label{beta_p} 
\end{align}
where ${{R}_{p}}\ge 0$. By re-exponentiation, 
\[\beta p\left( k \right)\le \exp\left( -\frac{k(k-1)}{2A} \right)=\exp\left( -
\frac{{{k}^{2}}}{2A}+\frac{k}{2A} \right)=\exp\left( -\frac{{{k}^{2}}}{2A} 
\right)\left[ 1+O\left( \frac{k}{a} \right) \right]
\]
so
\[p\left( k \right)\le \exp\left( -\frac{{{k}^{2}}}{2A} \right)\left[ 1+O\left( \frac{k}{a} \right) 
\right]\]
On the other hand, for $1\le k\le A/2$, by comparison with a geometric series,
\begin{align*}
 {{R}_{p}} &=\underset{j=0}{\overset{k-1}{\mathop 
\sum }}\,\left( \frac{{{j}^{2}}}{2{{A}^{2}}}+\frac{{{j}^{3}}}{3{{A}^{3}}}+\frac{{{j}^{4}}}{4{{A}^{4}}}\cdots  \right) \\ 
 & \le \underset{j=0}{\overset{k-1}{\mathop 
\sum }}\,\left( \frac{{{j}^{2}}}{2{{A}^{2}}}+\frac{{{j}^{3}}}{2{{A}^{3}}}+\frac{{{j}^{4}}}{2{{A}^{4}}}\cdots  \right) \\ 
 & =\underset{j=0}{\overset{k-1}{\mathop \sum }}\,\left( \frac{{{j}^{2}}}{2{{A}^{2}}}\frac{1}{1-j/A} 
\right)\le \underset{j=0}{\overset{k-1}{\mathop \sum }}\,\frac{{{j}^{2}}}{{{A}^{2}}}\le 
\frac{{{k}^{3}}}{{{A}^{2}}}  
\end{align*}
so by~\eqref{beta_p},
	\[\beta p\left( k \right)\ge \exp\left( -\frac{k(k-1)}{2A}-\frac{{{k}^{3}}}{{{A}^{2}}} \right)\ge 
\exp\left( -\frac{{{k}^{2}}}{2A}-\frac{{{k}^{3}}}{{{A}^{2}}} \right)\]
By~\eqref{nuisance_factor}
\begin{align*}
p\left( k \right)& \ge \exp\left( -\frac{{{k}^{2}}}{2A}-\frac{{{k}^{3}}}{{{A}^{2}}} 
\right)\left[ 1+O\left( \frac{k}{a} \right) \right] \\ 
 & =\exp\left( -\frac{{{k}^{2}}}{2A} \right)\exp\left( -\frac{{{k}^{3}}}{{{A}^{2}}} 
\right)\left[ 1+O\left( \frac{k}{a} \right) \right] \\ 
 & =\exp\left( -\frac{{{k}^{2}}}{2A} \right)~\left[ 1+O\left( \frac{k}{a}+\frac{{{k}^{3}}}{{{a}^{2}}} 
\right) \right] 
\end{align*}
\myqed
\end{proof}
The following lemma is similar to Gautschi's inequality \cite{Laforgia.1984fif}, but while the inequality
is restricted to $-1\le \nu \le 0$, the lemma here needs to hold for arbitrarily large negative $\nu$.
\begin{lemma}\label{factor_q}
The factor
	\[q\left( k \right)\triangleq \frac{\Gamma \left( k-\nu  \right)}{k!}\]
for $1\le k\le A$ and $\nu \le 0$ satisfies
	\[q\left( k \right)={{k}^{-\nu -1}}\left[ 1+O\left( \frac{1}{k} \right) \right]\]
\end{lemma}
\begin{proof} 
By Stirling's approximation for $k>0$
\cite[§6.1.37-38]{abramowitz.stegun.1972hom},
\cite[Eq. 5.11.3]{olver.et.al.2010nho}
\[\Gamma \left( k \right)=\sqrt{\frac{2\pi }{k}}{{\left( \frac{k}{e} 
\right)}^{k}}\left[ 1+O\left( \frac{1}{k} \right) \right]\]
and the relation
	\[{{\left( 1-\frac{\nu }{k} \right)}^{k}}=\exp \left[ k~\ln\left( 1-\frac{\nu }{k} \right) \right]=\exp 
\left[ -\nu +O\left( \frac{1}{k} \right) \right]={{e}^{-\nu }}\left[ 1+O\left( \frac{1}{k} \right) \right]\]
 gives 
\begin{align*}
 \frac{\Gamma \left( k-\nu  \right)}{\Gamma \left( k 
\right)}& =\sqrt{\frac{k}{k-\nu }}{{\left( \frac{k-\nu }{e} \right)}^{k-\nu }}{{\left( \frac{e}{k} 
\right)}^{k}}\left[ 1+O\left( \frac{1}{k} \right) \right] \\ 
 & ={{\left( 1-\frac{\nu }{k} \right)}^{-\frac{1}{2}-\nu }}{{\left( 1-\frac{\nu }{k} 
\right)}^{k}}{{e}^{\nu }}{{k}^{-\nu }}\left[ 1+O\left( \frac{1}{k} \right) \right] \\ 
 & =\left[ 1+O\left( \frac{1}{k} \right) \right]{{e}^{-\nu }}{{e}^{\nu }}{{k}^{-
\nu }}\left[ 1+O\left( \frac{1}{k} \right) \right] \\ 
 & ={{k}^{-\nu }}\left[1+O\left( \frac{1}{k} \right) \right]
\end{align*}
so from the definition of $q(k)$,
	\[q\left( k \right)=\frac{\Gamma \left( k-\nu  
\right)}{k!}=\frac{\Gamma \left( k-\nu  
\right)}{k\,\Gamma \left( k \right)}={{k}^{-\nu -
1}}\left[ 1+O\left( \frac{1}{k} \right) \right]
\]
\myqed
\end{proof}

\begin{figure}
\includegraphics[clip, trim=1.5cm 0.5cm 2.5cm 15cm ]{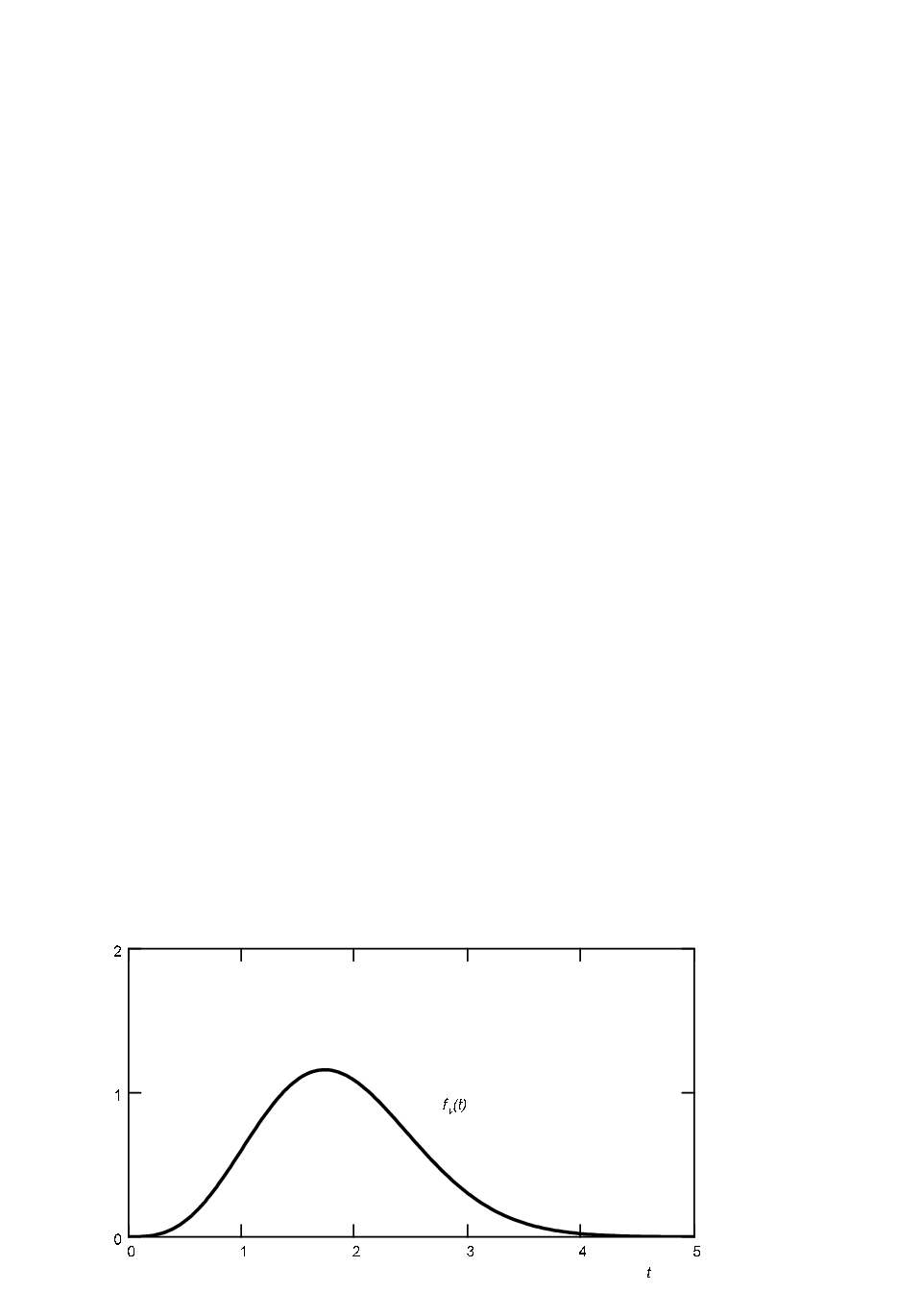}
\caption{The function ${{f}_{\nu }}(t)$ for $\nu =-4$.}
\label{fig:2}       
\end{figure}

Now it is time to take on the sum~\eqref{y_sum}, split in a head and tail part at index $M\triangleq 
\lceil{{A}^{3/4}}\rceil$, 
\begin{equation}\label{head_tail} 
\underset{k=0}{\overset{A}{\mathop \sum }}\,{{T}_{k}}=\underset{k=0}{\overset{M-
1}{\mathop \sum }}\,{{T}_{k}}+\underset{k=M}{\overset{A}{\mathop 
\sum }}\,{{T}_{k}}\triangleq {{R}_{head}}+{{R}_{tail}}
\end{equation}
Define $\Delta t\triangleq 1/\sqrt{A}$ and the function
\[{{f}_{\nu }}\left( t \right)\triangleq {{t}^{-\nu -1}}\exp\left( -\frac{{{t}^{2}}}{2} \right)\]
Clearly, the functions ${{f}_{\nu }}\left( t \right)$  (fig. \ref{fig:2})
and 
\[
f''_\nu(t) = \left[t^4 + (1+2\nu)t^2+(\nu^2+3\nu+2)\right] t^{-\nu-3} e^{-t^2/2}
\]
are continuous and bounded for bounded $\nu \le -3$ and $t\ge 0$.

\begin{lemma}\label{sum-approx}
The following relations hold for $\nu \le -3$:
\begin{equation}\label{sum-approx1}
\underset{k=M}{\overset{A~}{\mathop \sum }}\,{{f}_{\nu }}\left( k\Delta t 
\right)\Delta t=O\left(\frac{1}{\sqrt{a}} \right)
\end{equation}
and
\begin{equation}\label{sum-approx2}
\sum_{k=0}^A f_{\nu } \left( k\Delta t \right)\Delta t=
2^{-\nu /2-1} \Gamma \left( -\frac{\nu }{2} \right) 
+O\left(\frac{1}{\sqrt{a}} \right)
\end{equation}
\end{lemma}

\begin{proof}
According to the well-known trapezoidal rule, 
since ${{f}_{\nu }}\left( t \right)$ and ${f''_{\nu }}(t)$ are bounded 
for $\nu \le -3$ and $t\ge 0$, and for some $\tau \in 
\left[ M\Delta t,A\Delta t \right]$, 
\begin{align*}
\underset{M\Delta t}{\overset{A\Delta t}{\mathop 
\int }}\,{{f}_{\nu }}\left( t \right)dt 
& =
\left( \frac{f_\nu( A\Delta t)}{2}+\frac{f_\nu( M\Delta t)}{2} \right)\Delta t+
\\ & + \underset{k=M+1}{\overset{A-1~}{\mathop 
\sum }}\,{{f}_{\nu }}\left( k\Delta t \right)\Delta t-
\frac{A\Delta t-M\Delta t}{12}\Delta 
{{t}^{2}}{f''_\nu}\left( \tau  \right) \\ 
 & = 
\frac{f_\nu( A\Delta t)+f_\nu( M\Delta t)}{2}\,\Delta t
+\underset{k=M+1}{\overset{A-1~}{\mathop 
\sum }}\,{{f}_{\nu }}\left( k\Delta t \right)\Delta 
t+O\left( A \Delta t^3\right)  
\end{align*} 
so that
\[
\sum_{k=M}^A f_\nu\left( k\Delta t\right) \Delta t =
\int_{M\Delta t}^{A\Delta t}
f_\nu\left( t \right)dt 
+
O\left( \frac{1}{\sqrt{a}}\right)
\]
By substituting ${{t}^{2}}/2=u$ in the integral of ${{f}_{\nu }}$, the upper
incomplete gamma function
\cite[§6.5.3]{abramowitz.stegun.1972hom},
\cite[Eq. 8.2.2]{olver.et.al.2010nho}
 is obtained,
\begin{align}\label{integral1}
\mathop{\int }^{}{{f}_{\nu }}(t)dt & =\mathop{\int }^{}{{t}^{-\nu -1}}{{e}^{-\frac{{{t}^{2}}}{2}}}dt 
\notag \\ 
 & =\mathop{\int }^{}{{\left( \sqrt{2u} \right)}^{-\nu -2}}{{e}^{-u}}du \notag \\ 
 & ={{2}^{-\nu /2-1}}\mathop{\int }^{}{{u}^{-\nu /2-1}}{{e}^{-u}}du \notag \\ 
 & =-{{2}^{-\nu /2-1}}~\Gamma \left( -\frac{\nu }{2},u \right)+C  
\end{align}
Asymptotically
\cite[§6.5.32]{abramowitz.stegun.1972hom},
\cite[Eq. 8.11.2-3]{olver.et.al.2010nho},
\begin{equation}\label{asymp-gamma}
\Gamma \left( s,z \right)={{z}^{s-1}}{{e}^{-z}}\left[1+O\left( \frac{1}{z} 
\right) \right]
\end{equation}
implying that when $z$ increases, $\Gamma(s,z)$ approaches zero faster than any negative
power of $z$, including $1/\sqrt{a}$, i.e.,
\[
\Gamma \left( -\frac{\nu }{2},\frac{{{(M\Delta t)}^{2}}}{2} \right)
=
O\left( \frac{1}{\sqrt{a}}\right)
\]
This proves the first relation. For the second relation, by \eqref{integral1},
\begin{align*}
\int_0^{A\Delta t} f_\nu (t)dt 
&=
2^{-\nu /2-1} \Gamma \left( -\frac{\nu }{2}\right)
-
2^{-\nu /2-1} \Gamma \left( -\frac{\nu }{2},\frac{(A\Delta t)^2}{2} \right) \\
&=
2^{-\nu /2-1} \Gamma \left( -\frac{\nu }{2}\right)
+
O\left( \frac{1}{\sqrt{a}}\right)
\end{align*}
\myqed
\end{proof}

\begin{lemma}\label{Rtail} 
For bounded	 $\nu \le -3$, ${{R}_{tail}}=O\left( 1/\sqrt{a} \right)$.
\end{lemma}

\begin{proof} 
By Lemma~\ref{factor_p} and~\ref{factor_q},
\begin{align*}
0<~{{R}_{tail}}& = {{a}^{\nu /2}}\underset{k=M}{\overset{A~}{\mathop \sum }}\,q\left( k 
\right)p\left( k \right) \\ 
 & \le {{a}^{\nu /2}}\underset{k=M}{\overset{A~}{\mathop \sum }}\,{{k}^{-\nu -
1}}\left[1+ O\left( \frac{1}{k} \right) \right]~{{e}^{-{{k}^{2}}/2A}}\left[1+O\left( \frac{k}{a} \right) 
\right] \\ 
 & ={{a}^{\nu /2}}\underset{k=M}{\overset{A~}{\mathop \sum }}\,{{k}^{-\nu -1}}~{{e}^{-
{{k}^{2}}/2A}}O\left( 1 \right)  
\end{align*}
Substituting $k=k\Delta t \sqrt{A}$, 
\begin{align*}
 {{R}_{tail}}& ={{a}^{\nu /2}}\underset{k=M}{\overset{A~}{\mathop 
\sum }}\,{{\left( k\Delta t\sqrt{A} \right)}^{-\nu -1}}{{e}^{-
{{\left( k\Delta t \right)}^{2}}/2}}\Delta t\sqrt{A}\cdot 
O\left( 1 \right) \\ 
 & ={{\left( \frac{a}{A} \right)}^{\frac{\nu }{2}}}\underset{k=M}{\overset{A~}{\mathop 
\sum }}\,{{f}_{\nu }}\left( k\Delta t \right)\Delta t\cdot 
O\left( 1 \right) \\ 
 & =O\left( \frac{1}{\sqrt{a}} \right) 
\end{align*}
by Lemma~\ref{sum-approx}.
\myqed
\end{proof}

The term ${{R}_{head}}$
in~\eqref{head_tail} can be computed in a similar way. 

\begin{lemma}\label{Rhead} 
For bounded $\nu \le -4$,
\[
{{R}_{head}}={{2}^{-\nu /2-1}}\Gamma \left( -\frac{\nu }{2} 
\right)+O\left( \frac{1}{\sqrt{a}} \right)
\]
\end{lemma}

\begin{proof} 
This time $k<M$, and by Lemma~\ref{factor_p} and~\ref{factor_q},
\begin{align}
 {{R}_{head}}& ={{a}^{\nu /2}}\underset{k=0}{\overset{M-1}{\mathop \sum }}\,{{k}^{-\nu -
1}}\left[ 1+O\left( \frac{1}{k} \right) \right]{{e}^{-
\frac{{{k}^{2}}}{2A}}}\left[ 1+O\left( \frac{k}{a}+\frac{{{k}^{3}}}{{{a}^{2}}} \right) \right] \notag \\ 
 & ={{\left( \frac{a}{A} \right)}^{\nu /2}}\underset{k=0}{\overset{M-1}{\mathop 
\sum }}\,{{f}_{\nu }}\left( k\Delta t 
\right)\Delta t\left[ 1+O\left( \frac{1}{k}+\frac{k}{a}+\frac{{{k}^{3}}}{{{a}^{2
}}} \right) \right] \notag \\ 
 & =\underset{k=0}{\overset{M-1}{\mathop 
\sum }}\,{{f}_{\nu }}\left( k\Delta t 
\right)\Delta t\left[ 1+O\left( \frac{1}{a}+\frac{\Delta t
}{k\Delta t}+\frac{k\Delta t}{\sqrt{a}}+\frac{{{\left( k\Delta t \right)}^{3}}}{\sqrt{a}} \right) \right] \notag \\ 
 & =\underset{k=0}{\overset{M-1}{\mathop 
\sum }}\,{{f}_{\nu }}\left( k\Delta t 
\right)\Delta t+\underset{k=0}{\overset{M-1}{\mathop 
\sum }}\,{{f}_{\nu }}\left( k\Delta t \right)\Delta t\cdot 
O\left( \frac{\Delta t}{k\Delta t}+\frac{k\Delta t}{\sqrt{a}}+\frac{{{\left( k\Delta t \right)}^{3}}}{\sqrt{a}} 
\right) \notag \\ 
 & \triangleq S+\Delta S  \label{s_def}
\end{align}
Using the identity ${{f}_{\nu }}\left( t \right){{t}^{n}}={{f}_{\nu -n}}(t)$, the error term $\Delta{}S$ is
\begin{align*}
  \Delta S &=
\underset{k=0}{\overset{M-1}{\mathop \sum }}\,{{f}_{\nu 
+1}}\left( k\Delta t \right)\Delta t\cdot 
O\left( \Delta t \right)+ \\
&+
 \underset{k=0}{\overset{M-1}{\mathop 
\sum }}\,{{f}_{\nu -1}}\left( k\Delta t \right)\Delta t\cdot 
O\left( {1}/{\sqrt{a}} \right)+ \\ 
&+
 \underset{k=0}{\overset{M-1}{\mathop \sum }}\,{{f}_{\nu -
3}}\left( k\Delta t \right)\Delta t\cdot 
O\left({1}/\sqrt{a}\right) 
\end{align*}
Since
\[
0 \le \underset{k=0}{\overset{M-1}{\mathop \sum }}\,{{f}_{\nu 
}}\left( k\Delta t \right)\Delta t
\le
\underset{k=0}{\overset{A}{\mathop \sum }}\,{{f}_{\nu 
}}\left( k\Delta t \right)\Delta t
\]
which by Lemma~\ref{sum-approx} is bounded for $\nu \le -3$, 
\[
\Delta S = O\left( \frac{1}{\sqrt{a}} \right)
\]
for  $\nu \le -4$.
For the sum $S$ in~\eqref{s_def}, again using Lemma~\ref{sum-approx},
\begin{align*}
 S & =
\sum_{k=0}^{M-1}f_\nu\left( k\Delta t \right)\Delta t \\
&=
\sum_{k=0}^{A}f_\nu\left( k\Delta t \right)\Delta t -
\sum_{k=M}^{A}f_\nu\left( k\Delta t \right)\Delta t \\
 & ={{2}^{-\nu /2-1}}\Gamma \left( -\frac{\nu }{2} \right)+O\left( \frac{1}{\sqrt{a}}\right) 
\end{align*} \myqed
\end{proof}

By~\eqref{y_sum}, and combining Lemma~\ref{Rhead} and~\ref{Rtail},
\[
y_{\nu }^{a}\left( 0 \right)=\frac{{{2}^{\nu /2}}}{\Gamma \left( -\nu  
\right)}\left( {{R}_{head}}+{{R}_{tail}} \right)=\frac{\Gamma \left( -
\frac{\nu }{2} \right)}{2~\Gamma \left( -\nu  \right)}+O\left( \frac{1}{\sqrt{a}}\right)
\]
By the gamma function duplication rule
\cite[Eq. 1.2.3]{lebedev.1972sfa},
\cite[Eq. 5.5.5]{olver.et.al.2010nho},
\[
\frac{\Gamma \left( z \right)}{\Gamma \left( 2z 
\right)}=\frac{{{2}^{1-
2z}}~\sqrt{\pi }\text{ }\!\!~ }{\Gamma \left( z+\frac{1}{2} \right)}
\]
substituting $z=-\nu /2$,
\begin{equation}\label{y_zero} 
y_{\nu }^{a}\left( 0 \right)=
\frac{{{2}^{\nu }}\sqrt{\pi }}{\Gamma\left( \frac{1-
\nu }{2} \right)}+O\left( \frac{1}{\sqrt{a}} 
\right)={{H}_{\nu }}(0)+O\left( \frac{1}{\sqrt{a}} \right)
\end{equation}

\subsection{Convergence for $x=0$ and arbitrary $\nu$}

\begin{lemma}\label{y_zero_case} 
For $\nu$ in any bounded interval,
\[
{{y}_{\nu }}\left( 0 \right)={{H}_{\nu }}\left( 0 \right)+O\left( \frac{1}{\sqrt{a}} \right)
\]
and for $\Delta x=1/\sqrt{2a}$,
\[
   \frac{{{y}_{\nu }}\left( 0 \right)-{{y}_{\nu }}\left( -\Delta x 
\right)}{\Delta x}={H'_\nu\left( 0 \right)}+O\left( \frac{1}{\sqrt{a}} 
\right)
\]
\end{lemma}

\begin{proof} 
Given that ${{y}_{\nu }}\left( 0 \right)={{H}_{\nu }}(0)+O\left( 1/\sqrt{a} 
\right)$ for bounded $\nu <{{\nu }_{0}}$, then for ${{\nu }_{0}}\le \nu <{{\nu }_{0}}+1$ the 
difference equation~\eqref{difference_equation} can be rewritten into
\begin{equation}\label{rewritten-difference-eq}
{{c}_{n}}\left( \nu +1 \right)=\frac{\nu +a-n}{a}{{c}_{n}}\left( \nu  \right)-
\frac{\nu }{a}{{c}_{n}}\left( \nu -1 \right)
\end{equation}
so that for $n=A=\lceil a \rceil$ and by~\eqref{y_zero}, 
\begin{align}
\label{induction1}
 {{\left( 2a \right)}^{(\nu +1)/2}}{{c}_{A}}\left( \nu +1 \right)& ={{\left( 2a \right)}^{(\nu 
+1)/2}}~\left[ \frac{\nu +a-A}{a}{{c}_{A}}\left( \nu  \right)-\frac{\nu }{a}{{c}_{A}}\left( \nu -1 
\right) \right] \notag \\ 
 & =~\sqrt{2a}\frac{\nu +a-A}{a}~{{y}_{\nu }}\left( 0 \right)-\frac{\nu }{a}2a{{y}_{\nu -1}}\left( 0 
\right) \notag \\ 
 & =O\left( \frac{1}{\sqrt{a}} \right)-2\nu {{y}_{\nu -1}}\left( 0 \right) \notag \\ 
 & =-2\nu H_\nu(0) +O\left( \frac{1}{\sqrt{a}} \right) \notag \\ 
 & ={{H}_{\nu +1}}(0)+O\left( \frac{1}{\sqrt{a}} \right)
\end{align}
By induction, ${{y}_{\nu }}\left( 0 \right)={{H}_{\nu }}\left( 0 \right)+O\left( {1}/{\sqrt{a}} \right)$ for $\nu$ in any bounded interval. Additionally, by the backward recurrence relation~\eqref{backward_recurrence_relation} and the derivative rule for the Hermite function~\eqref{derivative_rule},
\begin{align}
\label{induction2}
 \frac{{{y}_{\nu }}\left( 0 \right)-{{y}_{\nu }}\left( -\Delta x 
\right)}{\Delta x}& =\frac{{{\left( 2a \right)}^{\nu /2}}{{c}_{A}}\left( \nu  \right)-
{{\left( 2a \right)}^{\nu /2}}{{c}_{A+1}}\left( \nu  \right)}{1/\sqrt{2a}} \notag \\ 
 & =\frac{\nu }{a}{{\left( 2a \right)}^{\nu /2+1/2}}{{c}_{A}}\left( \nu -1 \right) \notag \\ 
 & =2\nu {{y}_{\nu -1}}\left( 0 \right) \notag \\ 
 & =2\nu {{H}_{\nu -1}}\left( 0 \right)+O\left( \frac{1}{\sqrt{a}} \right) \notag \\ 
 & ={H'_\nu}\left( 0 \right)+O\left( \frac{1}{\sqrt{a}} \right)  
\end{align} \myqed
\end{proof}

\subsection{Convergence for arbitrary $x$ and arbitrary $\nu$}

In order to prove that ${{y}_{\nu }}\left( x \right)$ in~\eqref{def_y} converges to the solution of the Hermite 
differential equation~\eqref{Hermite_equation} having initial conditions $y\left( 0 \right)={{H}_{\nu }}(0)$ 
and $y'(0)=H'_\nu(0)$, it can be rewritten in normal form as
\begin{equation}\label{normal_y_equation} 
\myvec{y}'={\mymat{A}}(x)\myvec{y}
\end{equation}
where $\myvec{y}(x)\triangleq {{\left( y\left( x \right),~{y}'(x) \right)}^{T}}$and
\[{\mymat{A}}\left( x 
\right)\triangleq \left( \begin{matrix}
   0 & 1  \\
   -2\nu  & 2x  \\
\end{matrix} \right)
\]
Let $r\triangleq \sqrt{2a}$, $\Delta x\triangleq 1/r$, 
and  ${{x}_{k}}\triangleq k\Delta x$. Define a Cauchy polygon $\myvec{u}(x)$ 
for the differential equation~\eqref{normal_y_equation} by linear interpolation between points 
$\left( {{x}_{k}},{{\myvec{u}}_{k}} \right)$, where ${{\myvec{u}}_{0}}=\myvec{y}(0)$ and
\begin{equation}
\label{Cauchy_polygon} 
   {{\myvec{u}}_{k+1}}\triangleq {{\myvec{u}}_{k}}+\Delta x~{\mymat{A}}({{x}_{k}})~{{\myvec{u}}_{k}}
\end{equation}

\begin{lemma}\label{Cauchy_pol_conv} 
For $x$ and $\nu $ in bounded intervals $[0,\xi ]$ and $[-\psi ,\psi ]$, respectively, the 
Cauchy polygon $\myvec{u}(x)$ converges uniformly to the Hermite function solution with an error 
bound 
\[
\left| \myvec{u}(x)-\myvec{y}(x) \right|\le O\left( \frac{1}{\sqrt{a}} \right)
\]
\end{lemma}

\begin{proof} 
The Euclidean norm $||{\mymat{A}}(x)||$ of ${\mymat{A}}$ in~\eqref{normal_y_equation} equals the largest singular
value of the matrix, so
\[
||{\mymat{A}}(x)||={{\sigma }_{max}}\left({\mymat{A}}(x) \right)\le \sqrt{\text{tr}\left( {\mymat{A}}{{\left( x \right)}^{T}}{
\mymat{A}}(x) 
\right)}=\sqrt{1+4{{\nu }^{2}}+4{{x}^{2}}}
\]
Given arbitrary $\xi ,\psi >0$ and $L\triangleq ~\sqrt{1+4{{\psi }^{2}}+4{{\xi }^{2}}}$, for 
$x$ in $\left[ 0,\xi  \right]$ and $\nu $ in $[-\psi ,\psi ]$, by the definition of the Euclidean 
norm,
\[
\frac{\left| {\mymat{A}}\left( x \right)\left( \myvec{y}-\myvec{z} \right) \right|}{\left| \myvec{y}-\myvec{z} \right|}\le|| {\mymat{A}}(x)||\le L
\]
so $L$ is also a Lipschitz constant for~\eqref{normal_y_equation} when $x\in \left[ 0,\xi  \right]$. 
A definition and two theorems proved in~\cite[Sect. 7.3]{birkhoff.rota.1969ode} are now handy:
\begin{definition}
A vector function $\myvec{u}(x)$ is an {\em approximate solution} with {\em deviation} at most $\epsilon$
in the interval $a \le x \le \xi+a$ of the vector differential equation
\[
d\myvec{y}/dx = \myvec{Y}(\myvec{y},x),{\hskip 1cm}    a\le x\le a+\xi
\]
when $\myvec{u}(x)$ is continuous and satisfies the differential inequality
\[
\left| \myvec{u}'(x)-\myvec{Y}(\myvec{u}(x),x)\right|\le \epsilon
\]
for all except a finite number of points $x$ of the interval $\left[a, a+\xi \right]$.
\end{definition}

\begin{theorem}[Birkhoff and Rota, Th. 7.1]\label{BR_1}
Let the continuously differentiable function $\myvec{Y}$ satisfy $\left| \myvec{Y}\right| \le M$,
$\left| \partial\myvec{Y}/\partial x\right| \le C$, and $L$ be a
Lipschitz constant in the cylinder $D: \left| \myvec{y}-\myvec{c}\right| \le K$, $a \le x \le a + \xi$.
Then any Cauchy polygon in $D$ with partition $\pi$ is an approximate solution of $\myvec{y}'(x)=\myvec{Y}(\myvec{y},x)$
with deviation at most $(C + LM)|\pi|$.
\end{theorem}

\begin{theorem}[Birkhoff and Rota, Th. 7.3]\label{BR_2}
Let $\myvec{y}(x)$ be an exact solution and $\myvec{u}(x)$ be an approximate solution,
with deviation $\epsilon$, of the differential equation $\myvec{y}'(x)=\myvec{Y}(\myvec{y},x)$.
Let $\myvec{Y}$ satisfy a Lipschitz condition with Lipschitz constant $L$. Then, for $x \ge a$,
\[
\left| \myvec{y}(x) - \myvec{u}(x) \right| \le \left| \myvec{y}(a) - \myvec{u}(a)\right| e^{L(x-a)}
+ \left( \epsilon / L \right) \left( e^{L(x-a)} - 1\right)
\]
\end{theorem}

Bounds for $\left| {\mymat{A}}\left( x \right)\myvec{u}(x) \right|$ and $\left| 
\partial \left( {\mymat{A}}\left( x \right)\myvec{u}(x) \right)/\partial x \right|$ in $[0,\xi]$ can be chosen
\begin{equation}\label{rhs_power}
\left| {\mymat{A}}\left( {{x}_{k}} \right){{\myvec{u}}_{k}} \right|\le L\left| {{\myvec{u}}_{k}} \right|\le L\left| {{\myvec{u}}_{0}} 
\right|\underset{j=0}{\overset{k-1}{\mathop 
\prod }}\,||I+\Delta x~{\mymat{A}}\left( {{x}_{j}} \right)||\le L\left| {{\myvec{u}}_{0}} 
\right|{{\left( 1+L\Delta x~ \right)}^{k}}\le L\left| {{\myvec{u}}_{0}} 
\right|{{e}^{L\xi }}\triangleq M
\end{equation}
and
\[
{{\left| \frac{\partial \left( {\mymat{A}}\left( x \right)\myvec{u}\left( x \right) \right)}{\partial x} 
\right|}_{x={{x}_{k}}}}=\left| \left( \begin{matrix}
   0 & 0  \\
   0 & 2  \\
\end{matrix} \right){{\myvec{u}}_{k}} \right|\le 2\left| {{\myvec{u}}_{k}} \right|\le 2\left| {{\myvec{u}}_{0}} 
\right|{{e}^{L\xi }}\triangleq C
\]
By Theorem~\ref{BR_1}, Theorem~\ref{BR_2}, and Lemma~\ref{y_zero_case}, for $x\in [0,\xi ]$ and $\nu \in [-\psi ,\psi ]$,
\begin{align}\label{Cauchy_polygon_error} 
\left| \myvec{u}(x)-\myvec{y}(x) \right| &\le \left| \myvec{u}(0)-\myvec{y}(0) 
\right|{{e}^{Lx}}+\Delta x\left( \frac{C}{L}+M \right)\left( {{e}^{Lx}}-1 
\right)  \notag \\
 &=O\left( \frac{1}{\sqrt{a}} \right)
~{{e}^{L\xi }}+\Delta x\left( \frac{2}{L}+L \right)\left| 
{{\myvec{u}}_{0}} \right|{{e}^{L\xi }}\left( {{e}^{L\xi }}-1 \right) \notag \\
&=O\left( \frac{1}{\sqrt{a}} \right)
\end{align}
which is independent of $x$ and $\nu$,
so the Cauchy polygon~\eqref{Cauchy_polygon} converges uniformly to the 
Hermite function when $a\to \infty$. \myqed
\end{proof}
\par
Define $\myvec{z}_0\triangleq \myvec{u}_0$ and
\begin{align}
\label{z-def}
   {\myvec{z}_{k+1}} & \triangleq \left( {
\begin{matrix}
   {{y}_{\nu }}\left( {{x}_{k+1}} \right)  \\
   \dfrac{{{y}_{\nu }}\left( {{x}_{k+1}} \right)-{{y}_{\nu }}\left( {{x}_{k}} 
\right)}{\Delta x}  \\
\end{matrix} }
\right) \\
& = {{\myvec{z}}_{k}}+\Delta x\left( {\begin{matrix}
   \dfrac{{{y}_{\nu }}\left( {{x}_{k+1}} \right)-{{y}_{\nu }}\left( {{x}_{k}} 
\right)}{\Delta x}  \\
   \dfrac{{{y}_{\nu }}\left( {{x}_{k+1}} \right)-2{{y}_{\nu }}\left( {{x}_{k}} 
\right)+{{y}_{\nu }}\left( {{x}_{k-1}} \right)}{{\Delta{x}^{2}}} \notag \\
\end{matrix}}
\right)
\end{align}

Let $m\triangleq \lceil a-{{x}_{k}}r\rceil=a-{{x}_{k}}r+\left(\lceil a \rceil -a \right)=a-{{x}_{k}}r+\theta$, 
where $0\le \theta <1$. For simplicity of notation, the argument of ${{c}_{m}}$ is dropped when it is $\nu$. 
Consequently, 
\[{{\myvec {z}}_{k}}=
\left( \begin{matrix}
{{r}^{\nu }}{{c}_{m}} \\
{{r}^{\nu +1}}\left( {{c}_{m-1}}-{{c}_{m}} \right) \\

\end{matrix} \right)
\]
and
\[
   {{\myvec{z}}_{k+1}}={{\myvec{z}}_{k}}+\Delta x\left({\begin{matrix}
   {{r}^{\nu +1}}\left( {{c}_{m-1}}-{{c}_{m}} \right)  \\
   {{r}^{\nu +2}}\left( {{c}_{m-1}}-2{{c}_{m}}+{{c}_{m+1}} \right)  \\
\end{matrix}} \right)
\]
Multiplying the three-term recurrence relation~\eqref{three_term_recurrence_rule} by two, and substituting $x=\nu$ and $m=n$ 
gives the identity
\[
-2\nu {{c}_{m}}=2a{{c}_{m+1}}-\left( 2m+2a \right){{c}_{m}}+2m{{c}_{m-1}}\]
Rearranging, and using the facts that $2a={{r}^{2}}$ and $m=a-{{x}_{k}}r+\theta $,
\begin{equation}
\label{ttr-rearr}
{{r}^{2}}{{c}_{m-1}}-2{{r}^{2}}{{c}_{m}}+{{r}^{2}}{{c}_{m+1}}=2{{x}_{k}}r\left( {{c}_{m-1}}-
{{c}_{m}} \right)-2\nu {{c}_{m}}-2\theta \left( {{c}_{m-1}}-{{c}_{m}} \right)
\end{equation}
by which
\begin{align*}
  {{\myvec{z}}_{k+1}}& ={{\myvec{z}}_{k}}+\Delta x\left( \begin{matrix}
   {{r}^{\nu +1}}\left( {{c}_{m-1}}-{{c}_{m}} \right)  \\
   2\left( {{x}_{k}}-{\theta }/{r} \right){{r}^{\nu +1}}\left( {{c}_{m-1}}-{{c}_{m}} \right)-2\nu 
{{r}^{\nu }}{{c}_{m}}  \\
\end{matrix} \right) \\ 
 & ={{\myvec{z}}_{k}}+\Delta x\left( \begin{matrix}
   0 & 1  \\
   -2\nu  & 2{{x}_{k}}-2\theta /r  \\
\end{matrix} \right){{\myvec{z}}_{k}} \\ 
 & ={{\myvec{z}}_{k}}+\Delta x~{\mymat{A}}\left( {{x}_{k}} \right){{\myvec{z}}_{k}}+\Delta x\,\left( \begin{matrix}
   0 & 0  \\
   0 & -2\theta /r  \\
\end{matrix} \right){{\myvec{z}}_{k}}  
\end{align*}
This is nearly the same expression as for the Cauchy polygon~\eqref{Cauchy_polygon}, with only the $\theta$-term 
differing. Understanding the product sign below to 
multiply matrices in the proper order, and $\bm{I}$ to denote the identity matrix,
\begin{align*}
\frac{\left| {{\myvec{z}}_{k+1}}-{{\myvec{u}}_{k+1}} \right|}{\left| {{\myvec{u}}_{0}} \right|}
& \le 
\left\|
\underset{j=0}{\overset{k}{\mathop 
\prod }}\left[\mymat{I}+\Delta x~{\mymat{A}}\left( {{x}_{j}} 
\right)+\Delta x\left( \begin{matrix}
   0 & 0  \\
   0 & -2\theta /r  \\
\end{matrix} \right) \right]-\underset{j=0}{\overset{k}{\mathop 
\prod }}\left[\mymat{I}+\Delta x\,{\mymat{A}}\left( {{x}_{j}} \right) \right] 
\right\|
\end{align*}
Bounding the factor $\left\|{\mymat{I} + \Delta x~{\mymat{A}}{\left( x_j \right)} }\right\| \le \exp(L\xi)$ in
the same way as in~\eqref{rhs_power},
\begin{align}
\frac{\left| {{\myvec{z}}_{k+1}}-{{\myvec{u}}_{k+1}} \right|}{\left| {{\myvec{u}}_{0}} \right|}
& \le 
 \underset{j=1}{\overset{k}{\mathop \sum }}\,\left( \begin{matrix}
   k  \\
   j  \\
\end{matrix} \right)\Delta {{x}^{j}}
{\left\| {
{\left( \begin{matrix}   0 & 0  \\   0 & -2\theta /r  \\ \end{matrix} \right)}
}\right\| ^j e^{L\xi} } \notag
\\ 
\label{cauchy-diff1}
& = \underset{j=1}{\overset{k}{\mathop \sum }}\,\left( \begin{matrix}
   k  \\
   j  \\
\end{matrix} \right){{\left( \frac{2\theta \Delta x}{r} \right)}^{j}}{{e}^{L\xi }} 
=O\left( \frac{\xi }{r} \right){{e}^{L\xi }}=O\left( \frac{1}{\sqrt{a}} \right)  
\end{align}
demonstrates that $\myvec{z}$ converges uniformly to $\myvec{u}$ for $x \in\left[ 0,\xi  \right]$ and $\nu \in 
[-\psi ,\psi ]$. The proof for the descending direction from $x=0$ is omitted, since it is exactly 
analogous. By \eqref{Cauchy_polygon_error} and Lemma \ref{Cauchy_pol_conv},
\begin{equation}
\label{cauchy-diff2}
\left| \myvec{z}_k - \myvec{y}(x_k)\right| \le \left| \myvec{z}_k - \myvec{u}_k \right| +
\left| \myvec{u}(x_k) - \myvec{y}(x_k) \right| = O\left( \frac{1}{\sqrt{a}} \right) 
\end{equation}
so for $x_k \le x < x_{k+1}$,
\begin{align}
\label{cauchy-diff3}
\left| y^a_\nu(x)-H_\nu(x) \right| 
& \le 
\left| y^a_\nu(x_k)-H_\nu(x_k) \right| +  \left| y^a_\nu(x_k)-y^a_\nu(x_{k+1}) \right| \notag \\
& \le
\left| \myvec{z}_k-\myvec{y}(x_k) \right| +  \left| \myvec{z}_k-\myvec{z}_{k+1} \right| \notag \\
& = O\left( \frac{1}{\sqrt{a}} \right)
\end{align}
where the right hand side is independent of $x$ and $\nu$ for these parameters in any 
bounded interval.
\par
In order to demonstrate the sharpness of the bound, choose $\nu=2$, any real $x$, and arbitrarily 
large $a$ such that $n=a-x\sqrt{2a}$ is integer. Since $c_{2}^{a}\left( n \right)=1-\left( 1+2a 
\right)n/{{a}^{2}}+{{n}^{2}}/{{a}^{2}}$ and ${{H}_{2}}\left( x \right)=4{{x}^{2}}-2$,
\begin{equation}\label{sharp1}
y_{2}^{a}\left( x \right)-{{H}_{2}}\left( x \right)=\frac{2x\sqrt{2}}{\sqrt{a}}
\end{equation}
This completes the proof of Theorem \ref{theorem:1}. \qed

\section{Transition of the derivative}
\label{section:derivative}

\begin{theorem}
\label{theorem:conv-deriv}
For real $x$, $\nu$, and positive $a$,
\[
{\frac{\partial}{\partial\nu}
\left\{{
{\left( 2a \right)}^{\nu /2}}c_{\lceil a-x\sqrt{2a}\rceil}^{a}(\nu )
\right\}
=\frac{\partial}{\partial\nu}{{H}_{\nu }}\left( x 
\right)+O\left( \frac{1}{\sqrt{a}} \right)
}
\]
where $c_{n}^{a}(\nu )$ are Charlier polynomials and ${{H}_{\nu }}(x)$ is the Hermite 
function. The error bound $O\left( {1}/{\sqrt{a}}\; \right)$ is uniform for $\nu$ and $x$ in 
any bounded interval, and is sharp.
\end{theorem}

The proof of this theorem uses  same technique as the proof of Theorem \ref{theorem:1},
so the procedure can be abbreviated. First, the theorem is proved for the special case
$x=0$ and $\nu\le5$, then generalized to arbitrary $\nu$, and finally shown to converge
to the solution of a differential equation uniquely solved by the derivative of the Hermite function. 

\subsection{Convergence for $x=0$ and $\nu \le -5$}

Differentiating \eqref{y_sum} with respect to $\nu$,
\begin{equation}\label{T-diff1}
{\frac{\partial{y}_{\nu }(0)}{\partial\nu}}
=\frac{\partial}{\partial\nu} \left\{\frac{{{2}^{\nu /2}}}{\Gamma\left( -\nu  
\right)}\right\}\underset{k=0}{\overset{A}{\mathop \sum }}\,{{T}_{k}}
+
 \frac{{{2}^{\nu /2}}}{\Gamma\left( -\nu  
\right)}\underset{k=0}{\overset{A}{\mathop \sum }}\,\frac{\partial{{T}_{k}}}{\partial\nu}
\end{equation}
The first sum $\sum T_k$ is given by Lemmas \ref{Rtail} and \ref{Rhead}. Consider the second sum
\begin{equation}
\underset{k=0}{\overset{A}{\mathop \sum }}\,
\frac{\partial{{T}_{k}}}{\partial\nu}
=\underset{k=0}{\overset{M-
1}{\mathop \sum }}\,
\frac{\partial{{T}_{k}}}{\partial\nu}
+\underset{k=M}{\overset{A}{\mathop 
\sum }}\,
\frac{\partial{{T}_{k}}}{\partial\nu}
\triangleq {{R}_{head}}+{{R}_{tail}}
\end{equation}
Here
\begin{equation}\label{T-diff2}
\frac{\partial{{T}_{k}}}{\partial\nu}=
\left[ \ln \sqrt{a} - \psi(k - \nu)\right] T_k
\end{equation}
and $\psi(z)=\Gamma'(z)/\Gamma(z)$. By
\cite[§6.3.5 and §6.3.2]{abramowitz.stegun.1972hom},
\cite[Eqs. 5.4.12, 5.4.14, and 5.5.2]{olver.et.al.2010nho},
\[
\psi(k-\nu)
=
\psi(k) + O\left(\frac{1}{k}\right)
=
\ln{k}+O\left(\frac{1}{k}\right)
=\ln{k}\cdot\left[
1+O\left(\frac{1}{k \ln{k}}\right)
\right]
\]
Define
\[{{g}_{\nu }}\left( t \right)\triangleq \ln{t}\cdot f_\nu(t) = \ln{t}\cdot{{t}^{-\nu -1}}\exp\left( -\frac{{{t}^{2}}}{2} \right)
=- \frac{\partial}{\partial \nu} f_{\nu}(t)
\]
Since $t \ln t \to 0$ when $t \to 0^+$, taking zero as the value at $t = 0$,
the functions ${{g}_{\nu }}\left( t \right)$ and
\[
g''_\nu(t) =\ln t \cdot f''_\nu(t) - (3+2\nu+2t^2)f_{\nu-2}(t)
\]
are continuous and bounded for bounded $\nu \le -4$ and $t\ge 0$.

\begin{lemma}\label{Rtail2} 
For bounded	 $\nu \le -4$, ${{R}_{tail}}=O\left( 1/\sqrt{a} \right)$.
\end{lemma}

\begin{proof}
By Lemma~\ref{factor_p} and~\ref{factor_q}, and \eqref{T-diff2},
\begin{align*}
0<~{{R}_{tail}}
& = 
\frac{d}{d\nu}
\left\{
{{a}^{\nu /2}}\underset{k=M}{\overset{A~}{\mathop \sum }}\,
q\left( k 
\right)p\left( k \right) 
\right\}
\\ 
 & \le 
{{a}^{\nu /2}}\underset{k=M}{\overset{A~}{\mathop \sum }}\,
(\ln \sqrt{a} - \ln k)
\cdot {{k}^{-\nu -
1}}\left[1+ O\left( \frac{1}{k} \right) \right]~{{e}^{-{{k}^{2}}/2A}}\left[1+O\left( \frac{k}{a} \right) 
\right] \\ 
 & ={{a}^{\nu /2}}\underset{k=M}{\overset{A~}{\mathop \sum }}\,
(\ln \sqrt{a} - \ln k)
\cdot{{k}^{-\nu -1}}~{{e}^{-
{{k}^{2}}/2A}}O\left( 1 \right)  
\end{align*}
But substituting $k=k\Delta t \sqrt{A}$,
\begin{align*}
 {{R}_{tail}}
& =
{{a}^{\nu /2}}\underset{k=M}{\overset{A~}{\mathop 
\sum }}\,{
\left[\ln \sqrt{a} -\ln(k\Delta t\sqrt{A})\right]
 \cdot{\left( k\Delta t\sqrt{A} \right)}^{-\nu -1}}{{e}^{-
{{\left( k\Delta t \right)}^{2}}/2}}\Delta t\sqrt{A}\cdot 
O\left( 1 \right) \\ 
 & =
{{\left( \frac{a}{A} \right)}^{\frac{\nu }{2}}}\underset{k=M}{\overset{A~}{\mathop 
\sum }}\,
\left[\ln \sqrt{\frac{a}{A}} -\ln(k\Delta t)\right]
{f_{\nu }}
\left( k\Delta t \right)\Delta t\cdot 
O\left( 1 \right) \\ 
 & =
-\underset{k=M}{\overset{A~}{\mathop 
\sum }}\,
{g_{\nu }}
\left( k\Delta t \right)\Delta t\cdot 
O\left( 1 \right)
\end{align*}
Since $|\ln t| \le 1/t$ for $0 < t \le 1$ and $|\ln t| \le t$ for $t\ge 1$, $|g_\nu(t)| \le f_{\nu+1}(t) + f_{\nu-1}(t)$ for $t \ge 0$, and for $\nu \le -4$,
\[
|R_{tail}| \le 
\underset{k=M}{\overset{A~}{\mathop 
\sum }}\, {\left[f_{\nu+1 } \left( k\Delta t \right)+f_{\nu-1 }\left( k\Delta t \right)\right]} \Delta t\cdot O\left( 1 \right) = O\left( \frac{1}{\sqrt{a}} \right)
\]
by Lemma \ref{sum-approx}.
\myqed
\end{proof}

\begin{lemma}\label{Rhead2}
For bounded $\nu \le -4$,
\[
{{R}_{head}}=
\frac{d}{d\nu}
\left\{{{2}^{-\nu /2-1}}\Gamma \left( -\frac{\nu }{2} 
\right)\right\}+O\left( \frac{1}{\sqrt{a}} \right)
\]
\end{lemma}

\begin{proof}
This time $k<M$, and by Lemma~\ref{factor_p} and~\ref{factor_q}, like \eqref{s_def},
\begin{align}\label{new_s_def}
 {{R}_{head}}
& =
{{\left( \frac{a}{A} \right)}^{\nu /2}}\,
\sum_{k=0}^{M-1}\,
\left[\ln \sqrt{a} -\ln(k\Delta t\sqrt{A})\right]
\left( k\Delta t 
\right)\Delta t\left[ 1+O\left( \frac{1}{k}+\frac{k}{a}+\frac{{{k}^{3}}}{{{a}^{2
}}} \right) \right] \notag \\ 
 & =
-\underset{k=0}{\overset{M-1}{\mathop 
\sum }}\,{{g}_{\nu }}\left( k\Delta t 
\right)\Delta t
+
\underset{k=0}{\overset{M-1}{\mathop 
\sum }}\,{{g}_{\nu }}\left( k\Delta t \right)\Delta t\cdot 
O\left( \frac{\Delta t}{k\Delta t}+\frac{k\Delta t}{\sqrt{a}}+\frac{{{\left( k\Delta t \right)}^{3}}}{\sqrt{a}} 
\right) \notag \\ 
 & \triangleq S+\Delta S
\end{align}
Since $|{{g}_{\nu }}\left( t \right){{t}^{n}}|=|{{g}_{\nu -n}}(t)| \le f_{\nu - n +1}+f_{\nu-n-1}$, the error term $\Delta{}S$ is
\begin{align*}
 |\Delta S| &\le\underset{k=0}{\overset{M-1}{\mathop \sum }}
\,{\left[f_{\nu +2}\left( k\Delta t \right)+f_{\nu}\left( k\Delta t \right)\right]}\Delta t\cdot O\left( \Delta t \right)+ \\
&+ \underset{k=0}{\overset{M-1}{\mathop \sum }}
\,{\left[f_{\nu}\left( k\Delta t \right)+f_{\nu-2}\left( k\Delta t \right)\right]}\Delta t\cdot O\left( {1}/{\sqrt{a}} \right)+ \\ 
&+ \underset{k=0}{\overset{M-1}{\mathop \sum }}
\,{\left[f_{\nu -2}\left( k\Delta t \right)+f_{\nu -4}\left( k\Delta t \right)\right]}\Delta t\cdot O\left({1}/\sqrt{a}\right)=
O\left({1}/\sqrt{a}\right) 
\end{align*}
for $\nu \le -4$ by Lemma \ref{sum-approx}.

For the sum $S$ in~\eqref{new_s_def}, again using the trapezoidal rule,
as in the proof of Lemma~\ref{sum-approx},
\begin{align*}
 {{R}_{head}}
& =
-\underset{k=0}{\overset{M-1}{\mathop \sum }}\, g_\nu(t)
\left( k\Delta t \right) \Delta t \\
&=
-\int_0^{M\Delta t}g_\nu
\left(t \right) dt 
+ O\left( \frac{1}{\sqrt{a}} \right)  \\
&=
\int_0^{\infty}
\frac{d}{d\nu}
f_\nu
\left(t \right) dt 
+
\int_{M\Delta t}^{\infty}
g_\nu
\left(t \right) dt 
+ O\left( \frac{1}{\sqrt{a}} \right)  
\end{align*} \myqed
For $\nu \le -4$,
\[
\int_{M\Delta t}^{\infty}| g_\nu \left(t \right)| dt 
\le
\int_{M\Delta t}^{\infty} f_{\nu+1}\left(t \right) dt 
=
2^{-(\nu+3)/2} \Gamma\left(-\frac{\nu+1}{2},\frac{(M\Delta t)^2}{2}\right)
=
O\left( \frac{1}{\sqrt{a}} \right)  
\]
by \eqref{integral1} and \eqref{asymp-gamma}. 

The function $f_\nu$ satisfies $\int_0^\infty{f_{\nu} dt} < \infty$ for $\nu \le -1$, and
for $|h|\le 1$ and $t\ge 1$,
\[
\left(\frac{t^h - 1}{h}\right)
=
\left(\frac{e^{h \ln t} - 1}{h}\right)
=
\ln t + \frac{\ln^2 t}{2!}h + \frac{\ln^3 t}{3!}h^2+ \ldots
\le
e^{\ln t} - 1
< t
\]
The function $f_{\nu-1}(t)$ is an integrable function
dominating $\left|{f_{\nu+h}(t)-f_\nu(t)}\right|/{h}$ for $|h|\le 1$ and $t\ge 1$,
since
\[
\left|\frac{f_{\nu+h}(t)-f_\nu(t)}{h}\right| = \left(\frac{t^h - 1}{h}\right) f_\nu(t) < t f_{\nu}(t) = f_{\nu-1}(t)
\]
so by Lebesgue's dominant convergence theorem, the integration and differentiation
order can be switched in the integral
\[
\int_0^{\infty} \frac{d}{d\nu} f_\nu \left(t \right) dt
=
\frac{d}{d\nu} \left\{
\int_0^{\infty}  f_\nu \left(t \right) dt
\right\}
=
\frac{d}{d\nu} \left\{
2^{-\nu/2-1} \Gamma\left({-\frac{\nu}{2}}\right)
\right\}
\]
\myqed
\end{proof}

Using Lemma~\ref{Rhead2} and~\ref{Rtail2}, equation \eqref{T-diff1} becomes
\begin{align*}
\frac{d}{d\nu}
y_{\nu }^{a}\left( 0 \right)
&=
\frac{d}{d\nu}
\left\{ \frac{{2^{\nu /2}}}{\Gamma \left( -\nu  \right)} \right\}
\cdot
\sum_{k=0}^{A} T_k
+ \frac{{{2}^{\nu /2}}}{\Gamma \left( -\nu  \right)}
\left( {{R}_{head}}+{{R}_{tail}} \right) \\
&=
\frac{d}{d\nu}
\left\{ \frac{{2^{\nu /2}}}{\Gamma \left( -\nu  \right)} \right\}
\cdot
\left\{ 
2^{-\nu/2-1} \Gamma \left(-\frac{\nu}{2}\right)
+O\left( \frac{1}{\sqrt{a}} \right)
 \right\} \\
&+ \frac{{{2}^{\nu /2}}}{\Gamma \left( -\nu  \right)}
\cdot
\left\{
\frac{d}{d\nu}
\left[
2^{-\nu/2-1} \Gamma \left(-\frac{\nu}{2}  \right)
\right]
+ O\left( \frac{1}{\sqrt{a}}\right)
 \right\} \\
&=
\frac{d}{d\nu}
\left\{
\frac{\Gamma \left( -
\frac{\nu }{2} \right)}{2~\Gamma \left( -\nu  \right)}
\right\}
+O\left( \frac{1}{\sqrt{a}} \right) \\
&=
\frac{d}{d\nu}
H_\nu(0)+O\left( \frac{1}{\sqrt{a}} \right)
\end{align*}

\subsection{Convergence for $x=0$ and arbitrary $\nu$}
\begin{lemma}
For $\nu$ in any bounded interval,
\[
\frac{\partial{y}_{\nu }(0)}{\partial\nu}
=
{\frac{\partial{H}_{\nu }(0)}{\partial\nu}}+O\left( \frac{1}{\sqrt{a}} \right)
\]
and for $\Delta x=1/\sqrt{2a}$,
\[
\frac{\partial}{\partial\nu}
\left\{
\frac{{{y}_{\nu }(0)}-{y}_{\nu }(-\Delta x)}{\Delta x}
\right\}
=
\frac{\partial H'_\nu\left( 0 \right)}{\partial\nu}+O\left( \frac{1}{\sqrt{a}} \right)
\]
\end{lemma}

\begin{proof}
 Induction can be applied again, just as in the proof of Lemma \ref{y_zero_case}.
Given that ${\partial{y}_{\nu }(0)/\partial\nu}={\partial{H}_{\nu }}(0)/\partial\nu+O\left( 1/\sqrt{a} 
\right)$ for bounded $\nu <{{\nu }_{0}}$ and $n=A=\lceil a \rceil$,
again using \eqref{rewritten-difference-eq},
\begin{align*}
\frac{\partial}{\partial\nu}y_{\nu+1}(0)
&=
\frac{\partial}{\partial\nu}
\left\{
 {{\left( 2a \right)}^{(\nu +1)/2}}{{c}_{A}}\left( \nu +1 \right)
\right\} \\
& =
\frac{\partial}{\partial\nu}
\left\{
~\sqrt{2a}\frac{\nu +a-A}{a}~{{y}_{\nu }}\left( 0 \right)-\frac{\nu }{a}2a{{y}_{\nu -1}}\left( 0 \right) 
\right\} \\ 
& =
\frac{\sqrt{2}}{\sqrt{a}}
\frac{\partial}{\partial\nu}
\left\{
\nu y_{\nu}(0) \right\} +
\frac{(a-A)\sqrt{2}}{\sqrt{a}}\frac{\partial}{\partial\nu}
y_{\nu} \left( 0 \right) +
\frac{\partial}{\partial\nu}
\left\{
-2{\nu}{{y}_{\nu -1}}\left( 0 \right) 
\right\} \\ 
& =
O\left( \frac{1}{\sqrt{a}} \right) +
\frac{\partial}{\partial\nu}\left\{
-2\nu {{y}_{\nu -1}}\left( 0 \right) 
\right\}
\\ 
& =
-2 {{y}_{\nu -1}}\left( 0 \right) 
-2\nu \frac{\partial}{\partial\nu}\left\{
 {{y}_{\nu -1}}\left( 0 \right) 
\right\}
+ O\left( \frac{1}{\sqrt{a}} \right) 
\end{align*}
Applying the induction step,
\begin{align*}
\frac{\partial}{\partial\nu}y_{\nu+1}(0)
& =
-2 {{H}_{\nu -1}}\left( 0 \right) 
-2\nu \frac{\partial}{\partial\nu}\left\{
 {{H}_{\nu -1}}\left( 0 \right) 
\right\}
+ O\left( \frac{1}{\sqrt{a}} \right) \\
& =
\frac{\partial}{\partial\nu}\left\{
-2\nu {{H}_{\nu -1}}\left( 0 \right) 
\right\} + O\left( \frac{1}{\sqrt{a}} \right) \\
& =
\frac{\partial}{\partial\nu}\left\{
{{H}_{\nu +1}}\left( 0 \right) 
\right\} + O\left( \frac{1}{\sqrt{a}} \right)
\end{align*}
This implies that ${{y}_{\nu }}\left( 0 \right)={{H}_{\nu }}\left( 0 \right)+O\left( {1}/{\sqrt{a}} \right)$ for $\nu$ in any bounded interval. Using the
backward recurrence relation
\eqref{backward_recurrence_relation}, as in \eqref{induction2},
\begin{align*}
\frac{\partial}{\partial\nu}
\left\{
\frac{{{y}_{\nu }(0)}-{y}_{\nu }(-\Delta x)}{\Delta x}
\right\}
 & =
\frac{\partial}{\partial\nu}
2\nu {{y}_{\nu -1}}\left( 0 \right) \\ 
 & =
2 {{y}_{\nu -1}}\left( 0 \right) +
2\nu 
\frac{\partial}{\partial\nu}
{{y}_{\nu -1}}\left( 0 \right) \\ 
 & =
2 {{H}_{\nu -1}}\left( 0 \right) +
2\nu 
\frac{\partial}{\partial\nu}
{{H}_{\nu -1}}\left( 0 \right)+O\left( \frac{1}{\sqrt{a}} \right)  \\ 
 & =
\frac{\partial}{\partial\nu}\left\{
2\nu {{H}_{\nu -1}}\left( 0 \right) \right\}+O\left( \frac{1}{\sqrt{a}} \right) \\ 
 & =
\frac{\partial}{\partial\nu}
{H'_\nu}\left( 0 \right)+O\left( \frac{1}{\sqrt{a}} \right)  
\end{align*} \myqed
\end{proof}

\subsection{Convergence for arbitrary $x$ and arbitrary $\nu$}

By differentiating equation~\eqref{Hermite_equation} with respect to $\nu$, and defining $w \triangleq \partial y/\partial \nu$,
\begin{equation}
\label{diff-hermite}
w''=
2xw'-2\nu w-2y
=
2xw'-2\nu w-2H_\nu(x)
\end {equation}
This equation has the particular solution
$w(x) = \partial H_\nu(x)/\partial\nu$. The homogeneous equation is again the Hermite equation,
so the general solution of \eqref{diff-hermite} is
\[
w = \frac{\partial H_\nu(x)}{\partial\nu} + A H_\nu(x) + B H_\nu(-x)
\]
For initial conditions 
$w(0)=\partial H_\nu(0)/\partial\nu$ and $w'(0) = \partial H'_\nu(0)/\partial\nu$,
the unique solution of \eqref{diff-hermite} is obviously
$w(x) = \partial H_\nu(x)/\partial\nu$.

 Let $\myvec{y}(x)\triangleq {{\left( y\left( x \right),~{y}'(x),~w(x),~w'(x) \right)}^{T}}$.
Equation \eqref{diff-hermite} can be rewritten in the normal form \eqref{normal_y_equation}
where
\[{\mymat{A}}\left( x 
\right)\triangleq \left( \begin{matrix}
   0 & 1  & 0 & 0\\
   -2\nu  & 2x  & 0 & 0\\
  0 & 0 & 0 & 1  \\
 -2 & 0 &  -2\nu  & 2x  \\
\end{matrix} \right)
\]
This time, the Euclidean norm of ${\mymat{A}}$ satisfies 
$
||{\mymat{A}}(x)|| \le \sqrt{6+8{{\nu }^{2}}+8{{x}^{2}}}
$,
giving the Lipschitz constant $\sqrt{6+8{{\psi }^{2}}+8{{\xi }^{2}}}$.
In analogy with \eqref{Cauchy_polygon}, a Cauchy polygon $\myvec{u}(x)$ can be defined
such that for $x$ and $\nu $ in bounded intervals $[0,\xi ]$ and $[-\psi ,\psi ]$, respectively,
\[
{{\left| \frac{\partial \left( {\mymat{A}}\left( x \right)\myvec{u}\left( x \right) \right)}{\partial x} 
\right|}_{x={{x}_{k}}}}=\left| \left( \begin{matrix}
   0 & 0 & 0 & 0\\
   0 & 2 & 0 & 0 \\
   0 & 0 & 0 & 0\\
   0 & 0 & 0 & 2 \\
\end{matrix} \right){{\myvec{u}}_{k}} \right|\le 2\left| {{\myvec{u}}_{k}} \right|\le 2\left| {{\myvec{u}}_{0}} 
\right|{{e}^{L\xi }}
\]
As in the proof of Lemma \ref{Cauchy_pol_conv}, this means that
for bounded $x$ and $\nu$,
$\myvec{u}(x)$ converges uniformly to the solution $\myvec{y}$ with an error 
bound 
\[
\left| \myvec{u}(x)-\myvec{y}(x) \right|\le O\left( \frac{1}{\sqrt{a}} \right)
\]
Now, extending the definition of $\myvec{z}_k$ \eqref{z-def} to four components,
$\myvec{z}_0\triangleq \myvec{u}_0$ and
\[
   {\myvec{z}_{k+1}} \triangleq \left( {
\begin{matrix}
   {{y}_{\nu }}\left( {{x}_{k+1}} \right)  \\
   \dfrac{{{y}_{\nu }}\left( {{x}_{k+1}} \right)-{{y}_{\nu }}\left( {{x}_{k}} 
\right)}{\Delta x}  \\
   {{w}_{\nu }}\left( {{x}_{k+1}} \right)  \\
   \dfrac{{{w}_{\nu }}\left( {{x}_{k+1}} \right)-{{w}_{\nu }}\left( {{x}_{k}} 
\right)}{\Delta x}  \\
\end{matrix} }
\right) = {{\myvec{z}}_{k}}+\Delta x\left( {\begin{matrix}
   \dfrac{{{y}_{\nu }}\left( {{x}_{k+1}} \right)-{{y}_{\nu }}\left( {{x}_{k}} 
\right)}{\Delta x}  \\
   \dfrac{{{y}_{\nu }}\left( {{x}_{k+1}} \right)-2{{y}_{\nu }}\left( {{x}_{k}} 
\right)+{{y}_{\nu }}\left( {{x}_{k-1}} \right)}{{\Delta{x}^{2}}}  \\
   \dfrac{{{w}_{\nu }}\left( {{x}_{k+1}} \right)-{{w}_{\nu }}\left( {{x}_{k}} 
\right)}{\Delta x}  \\
   \dfrac{{{w}_{\nu }}\left( {{x}_{k+1}} \right)-2{{w}_{\nu }}\left( {{x}_{k}} 
\right)+{{w}_{\nu }}\left( {{x}_{k-1}} \right)}{{\Delta{x}^{2}}}  \\
\end{matrix}} \right)
\]
and writing $d_m \triangleq \partial c_m/\partial\nu$, 
\[{{\myvec {z}}_{k}}=
\left( \begin{matrix}
{{r}^{\nu }}{{c}_{m}} \\
{{r}^{\nu +1}}\left( {{c}_{m-1}}-{{c}_{m}} \right) \\
{{r}^{\nu }}{{d}_{m}} \\
{{r}^{\nu +1}}\left( {{d}_{m-1}}-{{d}_{m}} \right)
\end{matrix} \right)
\]
and
\[
   {{\myvec{z}}_{k+1}}={{\myvec{z}}_{k}}+\Delta x\left({\begin{matrix}
   {{r}^{\nu +1}}\left( {{c}_{m-1}}-{{c}_{m}} \right)  \\
   {{r}^{\nu +2}}\left( {{c}_{m-1}}-2{{c}_{m}}+{{c}_{m+1}} \right)  \\
   {{r}^{\nu +1}}\left( {{d}_{m-1}}-{{d}_{m}} \right)  \\
   {{r}^{\nu +2}}\left( {{d}_{m-1}}-2{{d}_{m}}+{{d}_{m+1}} \right)  \\
\end{matrix}} \right)
\]
Differentiating \eqref{ttr-rearr} with respect to $\nu$,
	\[{{r}^{2}}{{d}_{m+1}}-2{{r}^{2}}{{d}_{m}}+{{r}^{2}}{{d}_{m-1}}=2{{x}_{k}}r\left( {{d}_{m-1}}-
{{d}_{m}} \right)-2\nu {{d}_{m}}-2\theta \left( {{d}_{m-1}}-{{d}_{m}} \right)-2 c_m
\]
leads to
\begin{align*}
 {{\myvec{z}}_{k+1}}
& =
{{\myvec{z}}_{k}}+\Delta x\left( \begin{matrix}
   0 & 1  & 0 & 0\\
   -2\nu  & 2{{x}_{k}}-2\theta /r & 0 & 0 \\
  0 & 0 & 0 & 1 \\
  -2 & 0 & -2\nu  & 2{{x}_{k}}-2\theta /r  \\
\end{matrix} \right){{\myvec{z}}_{k}} \\ 
 & =
{{\myvec{z}}_{k}}+\Delta x~{\mymat{A}}\left( {{x}_{k}} \right){{\myvec{z}}_{k}}+
\Delta x
\left( \begin{matrix}
   0 & 0  & 0 & 0 \\
   0 & -2\theta /r & 0 & 0 \\
   0 & 0  & 0 & 0 \\
   0 & 0 & 0 & -2\theta /r \\
\end{matrix} \right){{\myvec{z}}_{k}}  
\end{align*}
By a procedure similar to the application of equations \eqref{cauchy-diff1}-\eqref{cauchy-diff3} in section \ref{section:conv},
\[
\left| \frac{\partial y^a_\nu(x)}{\partial\nu} - \frac{\partial H_\nu(x)}{\partial\nu} \right| \le O\left( \frac{1}{\sqrt{a}} \right)
\]
where the right hand side is independent of $x$ and $\nu$ for these parameters in any 
bounded interval.

The sharpness of the bound can be proved by contradiction:
Suppose that ${\partial y^a_{\nu}(x)}/{\partial\nu} -\partial H_\nu(x)/{\partial\nu} = O\left(b(a)\right)$
where $O(b(a))$ is tighter than $O\left(1/\sqrt{a}\right)$. Integrating this difference,
\[
\int_{\nu_1}^{\nu_2} \left[\frac{\partial y^a_{\nu}(x)}{\partial \nu} - {\frac{\partial H_\nu(x)}{\partial \nu}} \right]  d\nu
=
\left[y^a_{\nu_2}(x) - H_{\nu_2}(x)\right]
-
\left[y^a_{\nu_1}(x) - H_{\nu_2}(x)\right]
=O\left(b(a)\right)
\]
Choosing $\nu_1=0$ and $\nu_2=2$, arbitrary $x$, arbitrarily large $a$ such that
$a - x\sqrt{2a}$ is integer, and using \eqref{sharp1},
\[
\left[y^a_2(x) - H_2(x)\right]
-
\left[y^a_0(x) - H_0(x)\right]
=
\frac{2x\sqrt{2}}{\sqrt{a}} - (1 - 1)
= 
\frac{2x\sqrt{2}}{\sqrt{a}}
\]
which is a contradiction. This completes the proof of Theorem \ref{theorem:conv-deriv}.
\qed

\section{Convergence of zeros}
\label{section:zeros}

\begin{theorem}
For fixed real $x$ and positive $a\to \infty $, let $~n\triangleq ~\lceil a-x\sqrt{2a}\rceil$. 
For a convergent sequence of zeros ${{\nu }_{n}}\to \nu$ such 
that $c_{n}^{a}\left( {{\nu }_{n}} \right)=0$, the limit $\nu$ is a zero of the Hermite 
function, ${{H}_{\nu }}\left( x \right)=0$, satisfying $\nu ={{\nu }_{n}}+O\left( {1}/{\sqrt{a}}\; 
\right)$. Conversely, for a positive real zero $\nu$ of the Hermite function, there is a convergent 
sequence ${{\nu }_{n}}\to \nu$  of zeros of $c_{n}^{a}$ satisfying $\nu 
={{\nu }_{n}}+O\left( {1}/{\sqrt{a}}\; \right)$.
\end{theorem}

\begin{figure}
\includegraphics{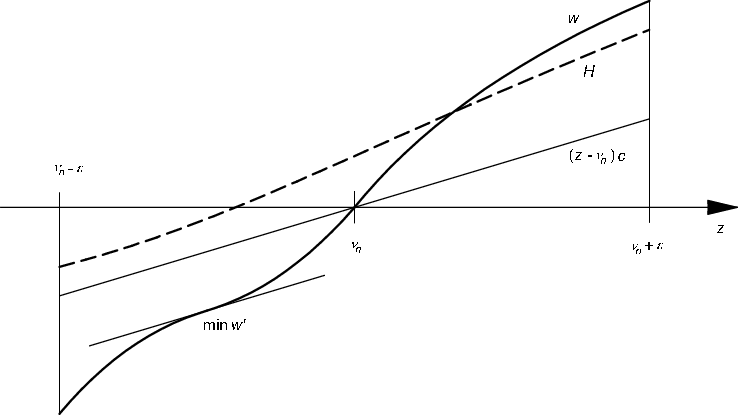}
\caption{The Hermite function must have a zero near the Charlier polynomial zero.}
\label{fig:3}       
\end{figure}

\begin{proof}
Define ${{w}_{n}}\left( z \right)\triangleq {{\left( \sqrt{2a} \right)}^{n}}c_{n}^{a}(z)$ and 
note that ${{w}_{n}}$ has the same zeros in $z$ as $c_{n}^{a}$. The proof is based on the 
well-known fact that the zeros of a Charlier polynomial are real, simple, and
positive~\cite{Kijima.1990otl}. Taylor-expanding ${{w}_{n}}\left( z \right)$ around one of its zeros $z={{\nu }_{n}}$, 
writing ${w'_n}({{\nu }_{n}})$ for $\partial {{w}_{n}}\left( z \right)/\partial z$ at $z={{\nu }_{n}}$,
\[
{{w}_{n}}\left( {{\nu }_{n}}+\varepsilon  \right)={{w}_{n}}\left( {{\nu }_{n}} \right)+\varepsilon 
{{w'_n}}({{\nu }_{n}})+O\left( {{\varepsilon }^{2}} \right)=\varepsilon 
\left( w'_{n}({{\nu }_{n}})\text{ }\!\!~ +O\left( \varepsilon  \right) \right)\triangleq 
\varepsilon W\left( {{\nu }_{n}},\varepsilon  \right)
\]
Since the zeros of a Charlier polynomial are simple, 
${{w'}_{n}}({{\nu }_{n}})\text{ }\!\!~ \ne 0$, the expression 
$W\left( {{\nu }_{n}},\varepsilon  \right)$ must be non-zero for $\varepsilon$ in some 
sufficiently small interval $I=[-\delta ,\delta ]$, where $0 < \delta \le {{\nu }_{n}}$. Assume that 
${w'_{n}}({{\nu }_{n}})>0$. The 
case $w'_{n}({{\nu }_{n}})\text{ }\!\!~ <0$ is treated in an analog way. 
Let $c\triangleq \underset{\varepsilon \in I}{\mathop{\inf }}\,W\left( {{\nu }_{n}},\varepsilon  
\right)$. Figure \ref{fig:3} illustrates $|(z-\nu_n)c|$ as a lower bound for $|w_n(z)|$.
By Theorem \ref{theorem:1}, due to the uniform convergence, for $z\in [{{\nu }_{n}}-
\delta ,{{\nu }_{n}}+\delta ]$, there is a $b$, independent of $n$ and $z$, such that
\[
\left| {{H}_{z}}\left( x \right)-{{w}_{n}}\left( z \right) \right|\le \frac{b}{\sqrt{a}}
\]
Choose $\varepsilon \triangleq \left( 1+b \right)/\left( c\sqrt{a} \right)$, which 
satisfies $\varepsilon <\delta$ for sufficiently large $a$. For
$z = {{\nu }_{n}}+\varepsilon $,
\[
{{H}_{z}}\left( x \right)\ge 
w_n\left( z\right)- \frac{b}{\sqrt{a}}
= \varepsilon W\left( {{\nu }_{n}},\varepsilon  \right)-
\frac{b}{\sqrt{a}}\ge \frac{1+b}{c\sqrt{a}}c-\frac{b}{\sqrt{a}}=\frac{1}{\sqrt{a}}>0
\]
Similarly, $z = {{\nu }_{n}}-\varepsilon$  implies that ${{H}_{z}}\left( x \right)<0$. 
Since ${{H}_{z}}(x)$ is an entire function and changes sign for z in $[{{\nu }_{n}}-
\varepsilon ,{{\nu }_{n}}+\varepsilon]$, it must have a zero there. By letting $a\to \infty$, the 
theorem is proved in one direction. For the reverse direction, switch the roles of $w$ and $H$. 
Assume that ${{H}_{\nu }}\left( x \right)=0$. Since ${{H}_{0}}(x)\equiv 1$, $\nu$ cannot be zero. 
Expand ${{H}_{z}}(x)$ around $z=\nu$, writing 
${\partial {{H}_{\nu}}(x)}/{\partial \nu}$ for $\partial {{H}_{z}}\left( x \right)/\partial z$ at $z=\nu$,
\[
{{H}_{\nu +\varepsilon }}\left( x \right)={{H}_{\nu }}\left( x \right)+\varepsilon 
{\partial {{H}_{\nu}}(x)}/{\partial \nu}+O\left( {{\varepsilon }^{2}} \right)=\varepsilon 
\left( {\partial {{H}_{\nu}}(x)}/{\partial \nu}+O\left( \varepsilon  \right) \right)\triangleq \varepsilon 
Z\left( \nu ,\varepsilon  \right)
\]
Let $x(\nu)$ be defined as the $p$th zero in $x$ of $H_\nu(x)=0$. 
It is known that $x(\nu)$ is a strictly monotonic function of $\nu$ for $\nu \ge 0$,
so $dx/d\nu \ne 0$ \cite{Elbert.Muldoon.1999iam}. Differentiating the equation by $\nu$,
\[
\frac{\partial {{H}_{\nu}}(x)}{\partial \nu} + \frac{\partial {{H}_{\nu}}(x)}{\partial x}\frac{dx}{d\nu} = 0
\]
so obviously, ${\partial {{H}_{\nu}}(x)}/{\partial \nu} = 0$ if and only if ${\partial {{H}_{\nu}}(x)}/{\partial x} = 0$. But if the latter derivative is zero, then ${{H}_{\nu -1}}\left( x \right)=0$  by the derivative rule~\eqref{derivative_rule}, and according to the three-term recurrence
for Hermite functions~\eqref{Hermite_three_term_recurrence}, all derivatives of ${{H}_{z}}(x)$ would be zero at $z=\nu$, 
entailing that $H$, being analytic, would be identically zero. In other words, all positive real zeros $\nu$ of $H_{\nu}(x)$ are simple.
\par
Consequently, ${\partial {{H}_{\nu}}(x)}/{\partial \nu} \ne 0$, and similarly to the first half of the proof,
$Z\left( \nu ,\varepsilon  \right)$ must be non-zero 
for $\varepsilon$ in some sufficiently small interval. It follows that ${{w}_{n}}\left( z 
\right)$ must be zero for some $z\in [\nu -\varepsilon ,\nu +\varepsilon ]$, where $\varepsilon 
=O\left( 1/\sqrt{a} \right)$. \myqed
\end{proof}

\section{Acknowledgements}
This research was funded by the European Union FP7 research project THE, 
"The Hand Embodied", under grant agreement 248587. The author is grateful for support by Dr. 
Henrik J\"orntell of Lund University, Dept. of Experimental Medical Science.





\bibliographystyle{model1b-num-names}










\end{document}